\documentclass{amsart}
\usepackage{amssymb, amsbsy, amsthm, amsmath, amstext, amsopn, verbatim}
\usepackage[all]{xy}
\usepackage{amsfonts}
\usepackage{amscd}

\newtheorem{thm}{Theorem} [section]
\newtheorem{lemma}[thm]{Lemma}
\newtheorem{cor}[thm]{Corollary}
\newtheorem{prop}[thm]{Proposition}
\newtheorem{notation}[thm]{Notation}

\theoremstyle{definition}

\newtheorem{defn}[thm]{Definition}

\newtheorem{example}[thm]{Example}

\theoremstyle{remark}

\newtheorem{remark}[thm]{Remark}

\newtheorem{convention}[thm]{Convention}

\begin{document}

\def\hangtwo{\hangindent=2.25em\hangafter=1}

\newcommand{\hs}{\mbox{\hspace{.4em}}}
\newcommand{\ds}{\displaystyle}
\newcommand{\bd}{\begin{displaymath}}
\newcommand{\ed}{\end{displaymath}}
\newcommand{\bcd}{\begin{CD}}
\newcommand{\ecd}{\end{CD}}

\newcommand{\proj}{\operatorname{Proj}}
\newcommand{\bproj}{\underline{\operatorname{Proj}}}
\newcommand{\spec}{\operatorname{Spec}}
\newcommand{\bspec}{\underline{\operatorname{Spec}}}
\newcommand{\pline}{{\mathbf P} ^1}
\newcommand{\pplane}{{\mathbf P}^2}

\newcommand{\ldb}{[[}
\newcommand{\rdb}{]]}

\newcommand{\Sym}{\operatorname{Sym}^{\bullet}}
\newcommand{\Symp}{\operatorname{Sym}}

\newcommand{\cA}{{\mathcal A}}
\newcommand{\bA}{{\mathbf A}}
\newcommand{\cB}{{\mathcal B}}
\newcommand{\cC}{{\mathcal C}}
\newcommand{\cs}{{\mathbf C} ^*}
\newcommand{\boldc}{{\mathbf C}}
\newcommand{\cD}{{\mathcal D}}
\newcommand{\cE}{{\mathcal E}}
\newcommand{\bE}{{\mathbf E}}
\newcommand{\cF}{{\mathcal F}}
\newcommand{\bF}{{\mathbf F}}
\newcommand{\cH}{{\mathcal H}}
\newcommand{\cI}{{\mathcal I}}
\newcommand{\cK}{{\mathcal K}}
\newcommand{\cL}{{\mathcal L}}
\newcommand{\bL}{{\mathbf L}}
\newcommand{\M}{{\mathcal M}}
\newcommand{\cM}{{\mathcal M}}
\newcommand{\bM}{{\mathbf M}}
\newcommand{\fM}{{\mathfrak m}}
\newcommand{\fm}{{\mathfrak m}}
\newcommand{\theo}{\mathcal{O}}
\newcommand{\boldp}{{\mathbf P}}
\newcommand{\boldq}{{\mathbf Q}}
\newcommand{\cQ}{{\mathcal Q}}
\newcommand{\bS}{{\mathbf S}}
\newcommand{\cU}{{\mathcal U}}
\newcommand{\bV}{{\mathbf V}}

\newcommand{\End}{\operatorname{End}}
\newcommand{\Hom}{\operatorname{Hom}}
\newcommand{\uHom}{\underline{\operatorname{Hom}}}
\newcommand{\Ext}{\operatorname{Ext}}
\newcommand{\bExt}{\operatorname{\bf{Ext}}}
\newcommand{\Tor}{\operatorname{Tor}}

\newcommand{\bk}{_{\bar{k}}}

\newcommand{\inv}{^{-1}}
\newcommand{\airtilde}{\widetilde{\hspace{.5em}}}
\newcommand{\airhat}{\widehat{\hspace{.5em}}}
\newcommand{\nt}{^{\circ}}
\newcommand{\perpen}{^{\bot}}

\newcommand{\supp}{\operatorname{supp}}
\newcommand{\id}{\operatorname{id}}
\newcommand{\res}{\operatorname{res}}
\newcommand{\lrar}{\leadsto}
\newcommand{\im}{\operatorname{Im}}

\newcommand{\Ker}{\operatorname{Ker}}

\newcommand{\TF}{\operatorname{TF}}
\newcommand{\Bun}{\operatorname{Bun}}
\newcommand{\Hilb}{\operatorname{Hilb}}
\newcommand{\nthord}{^{(n)}}
\newcommand{\Aut}{\underline{\operatorname{Aut}}}
\newcommand{\Gr}{\operatorname{\bf Gr}}
\newcommand{\Fr}{\operatorname{Fr}}
\newcommand{\pr}{$^\prime$}
\newcommand{\on}{\operatorname}
 
\def\hangfour{\hangindent=.4in\hangafter=1}
\def\hangup{\hangindent=3.2em\hangafter=1}
\def\hangtwo{\hangindent=.2in\hangafter=1}
\def\hangthree{\hangindent=.3in\hangafter=1}

\numberwithin{equation}{section}

\title[Descent of coherent sheaves]{Descent of coherent sheaves and complexes to geometric invariant theory quotients}
\author{Thomas Nevins}
\address{Department of Mathematics\\University of Illinois at Urbana-Champaign\\Urbana, IL 61801 USA}
\email{nevins@uiuc.edu}

\begin{abstract}
Fix a scheme $X$ over a field of characteristic 
zero that is equipped with an action of a reductive algebraic group $G$.  
   We give 
necessary and sufficient conditions for a $G$-equivariant coherent sheaf 
on $X$  or a bounded-above complex of $G$-equivariant coherent sheaves on 
$X$ to  descend to a good quotient $X/\!\!/G$.  This gives a 
description of the coherent 
derived category of $X/\!\!/G$ as an admissible subcategory of the equivariant
derived category of $X$.
\end{abstract}

\maketitle

\section{Introduction}
Varieties constructed using geometric invariant theory (or GIT) are
 ubiquitous in algebraic geometry.
 Many fundamental questions in the geometry of such GIT quotients concern the properties of certain natural vector bundles (or Chern classes of vector bundles) on them---for example, the classes 
$\kappa_i$ and $\lambda_j$   in the study of the geometry of $M_{g,n}$.
It is typical to construct such vector bundles on moduli spaces by 
{\em descent}.   
That is, one identifies the moduli space as a {\em good quotient}
(Definition \ref{goodquotientdef})
 $X/\!\!/G$ of a quasiprojective scheme
 $X$ by a reductive group $G$ using GIT.    Letting 
\bd
X^{ss}\xrightarrow{\pi} X/\!\!/G
\ed
denote the quotient map, 
one identifies a sheaf $\cM$ on the semistable locus $X^{ss}$ that one
expects to have the form $\cM = \pi^*\overline{\cM}$, where $\cM$ is the
desired sheaf on $X/\!\!/G$.  Verifying that such a $\overline{\cM}$ exists,
however, can be a nontrivial task: the morphism $\pi$ typically does not
satisfy hypotheses that would make it possible to apply Grothendieck's 
descent machinery.

 The following criterion, which may be found in \cite{MR90d:14008} (where the authors of that paper attribute it to Kempf), gives a convenient characterization of the vector bundles on $X^{ss}$ that descend to $X/\!\!/G$.

\begin{thm}\label{firsttheorem} \mbox{\rm (see \cite{MR90d:14008})}
Suppose $X$ is a quasiprojective scheme over an algebraically closed field $k$ of characteristic zero, and that $G$ is a reductive algebraic group over $k$ that acts on $X$ with a fixed choice of linearization $H$.  Let $E$ be a $G$-vector bundle on $X^{ss}$.  Then $E$ descends to $X/\!\!/G$ if and only if for every closed point $x$ of $X^{ss}$ such that the orbit $G\cdot x$ is closed in $X^{ss}$, the stabilizer of $x$ in $G$ acts trivially on the fiber $E_x$ of $E$ at $x$.
\end{thm}

\noindent
In this paper we extend the descent criterion of Theorem \ref{firsttheorem} in two directions.   First, we give a criterion for an arbitrary $G$-equivariant 
coherent sheaf to descend to a good quotient $X/\!\!/G$.
\begin{thm}\label{firstprop}
Suppose that $X$ is a scheme locally of finite type
 over a field $k$ of characteristic zero and that $G$ is a reductive algebraic group over $k$ that acts on $X$ with good quotient
$X\xrightarrow{\pi} X/\!\!/G$.    Let $\cM$ denote a $G$-equivariant coherent $\theo_X$-module.  Then the following are equivalent:

\vspace{.4em}

\noindent
\mbox{\rm (1)}\, $\cM$ descends to $X/\!\!/G$.

\hangthree\noindent
\mbox{\rm (2)}\, For every closed point $x\in X$ that lies in a closed
 $G$-orbit, the $\theo_{X,x}$-modules $\cM\otimes \theo_X/\fM_x$ and $\Tor_1^{\theo_X}\big(\cM,\theo_X/\fM_x\big)$ are generated by elements invariant under the isotropy subgroup $G_x$.

\vspace{.4em}

\noindent If $k$ is algebraically closed, these are equivalent also to

\vspace{.4em}

\noindent
\mbox{\rm (3)}\, For every closed point $x\in X$ that lies in a closed
 $G$-orbit, the $\theo_{X,x}$-modules $\cM\otimes\big(\theo_X/{\mathfrak m}_x\big)$ and $\Tor^{\theo_X}_1\big(\cM, \big(\theo_X/{\mathfrak m}_x\big)\big)$ are trivial representations of the isotropy group $G_x$.
\end{thm}

\noindent
Generalizing Theorem \ref{firsttheorem} in another direction, we also study the descent of equivariant complexes to $X/\!\!/G$.  Here
it makes more sense to consider the complex on $X$ as an object of the 
(equivariant) derived category,\footnote{The derived category we mean
 here is the derived category of quasicoherent
sheaves.  This need not coincide with the ``cohomologically quasicoherent''
derived category if $X$ is not quasi-compact and separated.}
 and to ask whether it descends up to 
quasi-isomorphism, i.e. whether it is in the essential image of the pullback
functor from the derived category of $X/\!\!/G$.
We then have the following theorem that describes the difference between the equivariant derived category of $X$---that is, the derived category of the quotient stack $[X/G]$---and the 
derived category of $X/\!\!/G$.

\begin{thm}\label{secondtheorem} 
Suppose $X$ is a scheme locally of finite type
 over a field $k$ of characteristic zero.  Suppose $G$ is a reductive algebraic group over $k$ that acts on $X$ with good quotient $X\xrightarrow{\pi} X/\!\!/G$.  
Let $\bE$
denote a bounded-above $G$-equivariant complex of coherent
 sheaves
 on $X$.  Then the following are equivalent:

\vspace{.4em}

\hangthree\noindent
\mbox{\rm (1)}\, ${\mathbf E}$ is equivariantly quasi-isomorphic to a complex ${\mathbf E'}$ on $X$ that descends to $X/\!\!/G$.

\hangthree\noindent
\mbox{\rm (2)}\, For every closed point $x\in X$ that lies in a closed 
$G$-orbit, the $\theo_X$-modules
\bd
H^j\big(\bE\overset{\bL}{\otimes}\theo_X/\fM_x\big)
\ed
are generated by elements invariant under the isotropy subgroup $G_x$ for all $j$.

\vspace{.4em}

\noindent
If $k$ is algebraically closed, these are equivalent also to

\vspace{.4em}

\hangthree\noindent
\mbox{\rm (3)}\,  For each 
 closed point $x\in X$ that lies in a closed $G$-orbit, the isotropy representations of $G_x$ on the $\theo_X$-modules
\bd
H^j\big(\bE\overset{\bL}{\otimes}\theo_X/\fM_x\big)
\ed
are trivial for all $j$.
\end{thm}
\noindent
In the case in which $G$ is finite, a similar result appeared in \cite{mathAG0206144}.

The author's original motivation for considering these descent problems was the possibility of applications in ``singular 
symplectic geometry.''  More precisely, important examples of singular symplectic moduli stacks come equipped with perfect pairings on
their tangent complexes that extend the symplectic structure (in an appropriate sense) on the smooth locus.  It is then natural to try to 
study the singularities of the associated coarse space in terms of the problem of descending this derived symplectic structure from the stack to the space.  We discuss some examples in this context in Section \ref{sympl section}.

It seems plausible that the descent results of this paper may be extended to characteristic
$p$, provided ``reductive'' is replaced by ``linearly reductive'' in appropriate places.  It would be nice to replace the hypothesis that
$G$ be linearly reductive in characteristic $p$ with a weaker hypothesis that stabilizers of closed points in closed orbits are linearly
reductive, but the author does not know how to prove, for example, Lemma \ref{canextendfromfiber} with such weakened hypotheses.

It is perhaps worth remarking that, for many interesting geometric
applications to GIT quotients $X/\!\!/G$, one needs to compare equivariant
sheaves on $X$ and $X^{ss}$, which seems to be extremely difficult 
 in general (see, for example, \cite{teleman}).

\vspace{.5em}

The author is grateful to
 Michel Brion, Mark Dickinson, Chris Skinner, and, especially, Igor
Dolgachev, for helpful discussions and comments.

\section{Preliminaries}

\subsection{Reminder about Group Actions}

\begin{convention}
Although it seems common  (see, for example, \cite{MR57:6035}) for connectedness of geometric fibers to be part of the definition of a reductive group scheme,
 in this paper we allow reductive groups to be disconnected; this causes us no trouble since we work always in characteristic zero.
\end{convention}

We remind the reader of a few facts we will need.

\begin{thm}\label{matsushimanisnevich} \mbox{\rm (Matsushima)}  Suppose $G$ is a reductive algebraic group over an algebraically closed field $k$ of characteristic $p\geq 0$, and $H$ is a closed $k$-subgroup of $G$.  Then the quotient scheme $G/H$ is affine if and only if
 (the identity component of) $H$ is reductive.
\end{thm}
\noindent
As a consequence, the stabilizers of closed points are reductive in our sense:
\begin{cor}\label{reductivebyfinite}
Suppose a reductive group $G$ over $k$ acts on an affine variety $X$ over $k$, and the closed point $x\in X$ has closed $G$-orbit in $X$.  Then the isotropy group $G_x$ is reductive.
\end{cor}

Recall also the following standard definition of a 
quotient arising in geometric invariant theory.

\begin{defn}\label{goodquotientdef}
Let $G$ be an affine algebraic group over $k$ acting on a $k$-scheme $X$.  A morphism $\phi: X\rightarrow Y$ is called a {\em good quotient}  if
\begin{itemize}
\item $\phi$ is affine and $G$-equivariant,
\item $\phi$ is surjective, and $U\subset Y$ is open if and only if $\phi^{-1}(U)\subset X$ is open,
\item the natural homomorphism $\theo_Y \rightarrow \big(\phi_*\theo_X\big)^G$ is an isomorphism,
\item if $W$ is an invariant closed subset of $X$, then $\phi(W)$ is a closed subset of $Y$; if $W_1$ and $W_2$ are disjoint invariant closed subsets of $X$, then $\phi(W_1)\cap\phi(W_2) = \emptyset$.
\end{itemize}
\end{defn}

\noindent
Note that any linearized action of a reductive group on a quasi-projective scheme $V$ over a field admits such a quotient of the semistable locus $X=V^{ss}$
(see Theorem 4 of \cite{MR57:6035}).

\subsection{Some Basic Facts}
Let
\bd
\pi: X \longrightarrow X/\!\!/G
\ed
denote the projection morphism.
Recall that a quasicoherent sheaf $\cM$ on $X$ {\em descends to $X/\!\!/G$}
if there exists a quasicoherent sheaf $\overline{\cM}$ on $X/\!\!/G$ and 
an isomorphism $\pi^*\overline{\cM}\rightarrow \cM$.  Similarly, a 
$G$-equivariant complex of quasicoherent sheaves $\bM$ on $X$ is said to
{\em descend to $X/\!\!/G$} if there is a complex $\overline{\bM}$ on
$X/\!\!/G$ and a quasiisomorphism ${\mathbf L}\pi^*\overline{\bM}\rightarrow
\bM$.

Recall that one defines the functor $\pi_*^G$ by $\pi_*^G(M) = (\pi_*M)^G$.
Since $\pi$ is affine and $G$ is reductive, $\pi_*^G$ is exact.
\begin{lemma}  The functors
$(\pi^*, \pi_*^G)$ form an adjoint pair,
\bd
\xymatrix{\on{qcoh}(X, G) \ar@<0.5ex>[r]^{\pi_*^G} & \on{qcoh}(X/\!\!/G),
 \ar@<0.5ex>[l]^{\pi^*}}
\ed
where $\on{qcoh}(X,G)$ denotes the category of $G$-equivariant
quasicoherent $\theo_X$-modules.  These induce an adjoint pair $({\mathbf L}\pi^*, \pi_*^G)$ of derived functors,
\bd
\xymatrix{D^{-}(\on{qcoh}(X,G)) \ar@<0.5ex>[r]^{\pi_*^G} & D^{-}(\on{qcoh}(X/\!\!/G)). \ar@<0.5ex>[l]^{{\mathbf L}\pi^*}}
\ed
\end{lemma}
  Furthermore, the functors $\pi^G_*\pi^*$ and $\pi^G_*{\mathbf L}\pi^*$ are
(isomorphic to) the identity functors.  We thus obtain:
\begin{cor}\label{canonical adjunction}
\mbox{}
\begin{enumerate}
\item A $G$-equivariant quasicoherent sheaf on $X$ descends to $X/\!\!/G$ 
if and only if the canonical map
\begin{equation}\label{invariantsgivedescent}
\pi^*\pi_*^G\cM \rightarrow \cM
\end{equation}
is an isomorphism.
\item A $G$-equivariant quasicoherent 
complex $\bM$ on $X$ descends to $X/\!\!/G$ if and only
if the canonical morphism 
\begin{equation}\label{derived cat canonical}
{\mathbf L}\pi^*\pi^G_*\bM\rightarrow \bM
\end{equation}
 is a 
quasi-isomorphism.
\end{enumerate}
\end{cor}
\begin{remark}\label{reduce to affine}
It is immediate from Corollary \ref{canonical adjunction} that the descent
criteria of Theorems \ref{firsttheorem} and \ref{secondtheorem} may be checked
locally on $X$.  Consequently, we may assume in the proofs that $X=\spec A$ is
affine.
\end{remark}
Let $\cK = \on{ker}(\pi_*^G)$ denote the kernel of the functor $\pi_*^G$: this
is the subcategory of objects $C$ for which $\pi_*^GC=0$.  Let $\cI$ denote the essential image of $\bL\pi^*$.  Standard arguments then give:

\begin{prop}
The pair $(\cK,\cI)$ forms a semiorthogonal decomposition of
$D^-(\on{qcoh}(X,G))$.
\end{prop}
Thus, Theorem \ref{secondtheorem} 
may be interpreted as a characterization of the
coherent part of the summand $\cI$ in this decomposition.

\begin{example}
Let $X=\spec {\mathbf C}[z]$, a variety over ${\mathbf C}$, with the usual
${\mathbf G}_m$-action.  Let $\cM = \widetilde{z^{-1}{\mathbf C}[z]}$.  
Then 
\bd
\pi^*\pi_*^G \cM = \widetilde{{\mathbf C}[z]}\subset 
\widetilde{z^{-1}{\mathbf C}[z]},
\ed
and the cone $C$
 on this map is the fiber of $\cM$ at the origin.  Thus, this gives
an example of a $C[-1]\in\on{ker}({\mathbf L}\pi^*\pi_*^G)$ and 
$\cF\in\on{im}({\mathbf L}\pi^*)$ such that 
$\on{Hom}(C[-1],\cF)\neq 0$ (in the derived category). So this semiorthogonal
decomposition is not an orthogonal decomposition. $\cM$ is also an
example of a sheaf whose fiber at every closed point $x$
 in a closed orbit contains {\em no} $G_x$-invariant elements, but 
$\pi_*^G\cM\neq 0$; so $\cK$ cannot be expected to have such a simple fiberwise
description as $\cI$.  
\end{example}

\begin{notation}
 We will denote by $\bar{k}$ a fixed algebraic closure of $k$, and by $X_{\bar{k}}$, $\cM_{\bar{k}}$ and so on the tensor products with $\bar{k}$. 
\end{notation}

\begin{lemma}\label{base-changed condition}
\mbox{}
\begin{enumerate}
\item
Suppose that $\cM$ satisfies hypothesis (2) of Theorem \ref{firstprop}.  Then
$\cM_{\bar{k}}$ satisfies hypothesis (3) of Theorem \ref{firstprop}.
\item Suppose that $\bM$ satisfies hypothesis (2) of Theorem 
\ref{secondtheorem}.  Then $\bM_{\bar{k}}$ satisfies hypothesis (3) of 
Theorem \ref{secondtheorem}.
\end{enumerate}
\end{lemma}
\begin{proof}
Let $x$ denote a closed point of $X_{\bar{k}}$ and let $y$ denote its image
in $X$; this is a closed point of $X$.
There is an inclusion of isotropy groups
\begin{equation}\label{isotropyinclusion}
\left(G_{\bar{k}}\right)_x \subseteq \left(G_y\right)_{\bar{k}};
\end{equation}
indeed, the $G_{\bar{k}}$-action on $X_{\bar{k}}$ covers the $G$-action on $X$,
 and so in particular the image of $\left(G_{\bar{k}}\right)_x$ in $G$ 
fixes $y$.  

(1) Since $\theo_{X_{\bar{k}}}/\fm_x=\bar{k}$ is flat over $\theo_X/\fm_y$, we have: 
for any $\theo_X$-module $N$ and any $i\geq 0$,
\begin{equation}\label{tors}
\Tor_i^{X_{\bar{k}}}\left(N_{\bar{k}},\theo_{X_{\bar{k}}}/\fm_x\right) =  \Tor_i^X\left(N,\theo_X/\fm_y\right)\underset{\theo_X/\fm_y}{\otimes} \bar{k}.
\end{equation}

Combining \eqref{isotropyinclusion} and \eqref{tors}, we find that $\cM_{\bar{k}}\otimes \theo_{X_{\bar{k}}}/\fm_x$ and $\Tor_1^{X_{\bar{k}}}\left(\cM_{\bar{k}}, \theo_{X_{\bar{k}}}/\fm_x\right)$ are generated over $\bar{k}$ by elements invariant under $\left(G_{\bar{k}}\right)_x$, and thus are in fact trivial representations of this isotropy group.

(2) Again, because $\theo_{X_{\bar{k}}}/\fm_x=\bar{k}$ is flat over $\theo_X/\fm_y$, we have
\bd
H^i(\bL\pi^*\bM\overset{\bL}{\otimes}\theo_X/\fm_x) = 
H^i(\bM\overset{\bL}{\otimes}\theo_{X/\!\!/G}/\fm_y\otimes\theo_X/\fm_x)
= H^i(\bM\overset{\bL}{\otimes}\theo_{X/\!\!/G}/\fm_y)\otimes\theo_X/\fm_x.
\ed
By hypothesis, $H^i(\bM\overset{\bL}{\otimes}\theo_{X/\!\!/G}/\fm_y)$ is generated by 
$G_y$-invariant elements.  Combined with \eqref{isotropyinclusion}, this implies that $H^i(\bL\pi^*\bM\overset{\bL}{\otimes}\theo_X/\fm_x)$ is generated
by $\theo_X/\fm_x$-invariants, as desired.
\end{proof}

\begin{lemma}\label{canextendfromfiber}
Suppose that $k$ is an algebraically closed field of characteristic zero
 and $A$ is a finitely 
generated $k$-algebra with a rational action of a reductive algebraic
$k$-group $G$.  Let $\cF$ be a $G$-equivariant 
coherent sheaf on $\spec A$.  Let $x\in \spec A$ be a closed point
 the $G$-orbit through which is closed in $\spec A$ and is defined by the ideal $I\subseteq A$.  If $s \in \cF\otimes \theo/\fm_x$ is $G_x$-invariant, then there is a $G$-invariant section $\tilde{s}\in \cF(\spec A)$ such that the image of $\tilde{s}$ in $\cF\otimes \theo/\fm_x$ is the element $s$.
\end{lemma}
\begin{proof}
Restricting $\cF$ to $\spec (A/I)$, we obtain a $G$-equivariant coherent
sheaf on $G/G_x$, which must therefore be the vector bundle associated to some representation $V$ of $G_x$.  Now $s$ determines an element of $V$, and by hypothesis $s$ lies in a trivial $G_x$-subrepresentation $W$ of $V$.  The associated bundle $G\times_{G_x} W$ is trivial, and thus $s$ extends to a $G$-invariant global section of the pullback of $\cF|_{\spec (A/I)}$.  But there is a $G$-equivariant surjection $\cF \rightarrow \cF|_{\spec(A/I)}$, and because $G$ is reductive, we may lift this section to a $G$-invariant section of $\cF$.
\end{proof}
\begin{remark}
The previous lemma easily extends to quasicoherent sheaves, but that extension
seems not to be as useful.
\end{remark}
\begin{lemma}\label{canfindsurjection}
Let $k$ be an algebraically closed field of characteristic zero.
Suppose $\cM$ is a $G$-equivariant quasicoherent sheaf on an affine
 $G$-scheme $Y$.  Then there is a $G$-equivariant locally free 
(quasicoherent) sheaf $\bV$ on $Y$ and a surjective $G$-equivariant homomorphism $\bV\rightarrow\cM$.  Moreover, if $\cM$ is coherent and
 the stabilizer $G_x$ acts trivially on the
fiber of $\cM$ at $x$ for every closed point
$x$ in a closed $G$-orbit, then $\bV$ may be chosen to be coherent and
to descend to $Y/\!\!/G$.
\end{lemma}
\begin{proof}
For the first statement we may take $\bV = \theo\otimes_k H^0(\cM)$.  
For the second statement, we start with
$\bV' =\theo\otimes_k H^0(\cM)^G$: the natural map to $\cM$ is $G$-equivariant,
and Lemma \ref{canextendfromfiber} implies that it is 
surjective on fibers at closed points in 
closed $G$-orbits.  Consequently, the cokernel is supported on a $G$-invariant
closed subset in the complement of the union of the closed $G$-orbits, 
implying 
that its support is empty.  Now, since $\cM$ is finitely generated, we 
may replace $H^0(\cM)^G$ by a finite-dimensional $k$-vector
 subspace $W$ and take
$\bV = \theo\otimes_k W$.
\end{proof}

\section{Proof of Theorem \ref{firstprop}}\label{proofofone}
Recall that, by Remark \ref{reduce to affine}, we may assume that $X$ is
affine.
We begin with:
\begin{lemma}\label{presentation descent}
A $G$-equivariant quasicoherent $\theo_X$-module $M$ descends to $X/\!\!/G$
if and only if there is a $G$-equivariant presentation
\bd
P_1\rightarrow P_0\rightarrow M\rightarrow 0,
\ed
where $P_0$ and $P_1$ are $G$-equivariant locally free sheaves
that descend to $X/\!\!/G$.
\end{lemma}
\begin{proof}
The ``only if'' part is immediate from right-exactness of  $\pi^*$ and 
the equation ${\mathbf 1}  = \pi_*^G\pi^*$ for the identity functor.  
Conversely, if $M$ has a presentation as 
above, then, letting $\overline{M} = \operatorname{coker}(\pi_*^G P_1
\rightarrow \pi_*^G P_0)$, \eqref{invariantsgivedescent} and right-exactness
of $\pi^*$ give $\pi^*\overline{M} = M$, as desired.
\end{proof}

\subsection{Necessity of the Criterion}
Suppose first that $\cM$ descends to $X/\!\!/G$; we will prove that condition (2) of Theorem \ref{firstprop} is satisfied.  
Let $x\in X$ be a closed point in a closed $G$-orbit, and 
let $G_x$ denote the isotropy subgroup of $x$.

Choose a presentation 
\bd
P_1\rightarrow P_0\rightarrow \pi_*^G\cM\rightarrow 0
\ed
of $\pi_*^G\cM$ by  locally free coherent sheaves on $X/\!\!/G$.  By 
Corollary \ref{canonical adjunction} and right-exactness of tensor product,
$\pi^*(P_1\rightarrow P_0)$ gives a presentation of $\cM$ by
locally free sheaves on $X$.  The fiber of $P_i$ at $\theo_X/\fm_x$ is given by
\begin{equation}\label{compute fiber}
\pi^*P_i\otimes\theo_X/\fm_x = (P_i\otimes\theo_{X/\!\!/G}/\fm_{\pi(x)})
\otimes \theo_X/\fm_x.
\end{equation}
These fibers are generated over $\theo_X/\fm_x$ by their subspaces 
$P_i\otimes\theo_{X/\!\!/G}/\fm_{\pi(x)}$, which consist of $G_x$-invariants.
Moreover, it follows from \eqref{compute fiber} that
\bd
\cK = \on{ker}((\pi^*P_1)\otimes\theo_X/\fm_x\rightarrow(\pi^*P_0)\otimes\theo_X/\fm_x)
\ed
 is defined over $\theo_{X/\!\!/G}/\fm_{\pi(x)}$; consequently, $\cK$
is also generated by $G_x$-invariants.
But the fiber
$\cM\otimes\theo_X/\fm_x$ is $G_x$-equivariantly isomorphic to a quotient
of the fiber of $\pi^*P_0$, and $\on{Tor}_1(\cM,\theo_X/\fm_x)$ is 
equivariantly isomorphic to a quotient of $\cK$.
So condition (2) is satisfied.

\subsection{Sufficiency of the Criterion}
Now, suppose that condition (2) holds.  We will show that this implies condition (1).

\vspace{.5em}

\noindent
{\bf Case 1.\,  $k$ algebraically closed.}\;
By assumption and Lemma \ref{canfindsurjection}, there is a surjective 
$G$-equivariant
homomorphism $P_0\rightarrow \cM$ where $P_0$ is a
locally free coherent sheaf that 
descends.
Let $\cK$ denote its kernel.
For every closed point 
$x\in X$ lying in a closed $G$-orbit, we have an exact sequence
\bd
\Tor_1(\cM,\theo_X/\fm_x) \rightarrow \cK\otimes \theo_X/\fm_x \rightarrow P_0 \otimes \theo_X/\fm_x \rightarrow 
\cM\otimes \theo_X/\fm_x \rightarrow 0.
\ed
Because $G_x$ acts trivially on $P_0\otimes\theo_X/\fm_x$ by construction and on $\Tor_1(\cM,\theo_X/\fm_x)$ by assumption,  $G_x$ acts trivially on $\cK\otimes\theo_X/\fm_x$ for each such $x$.  Applying Lemma \ref{canfindsurjection} to
 $\cK$, we get a $G$-equivariant homomorphism
\bd
P_1 \longrightarrow P_0,
\ed
the image of which is $\cK$, where $P_1$ is a $G$-equivariant locally free 
sheaf that descends.  Thus, we have a $G$-equivariant presentation
$P_1 \rightarrow P_0 \rightarrow \cM \rightarrow 0$ 
in which $P_1$ and $P_0$ descend to $X/\!\!/G$. 
By Lemma \ref{presentation descent}, this 
completes the proof when $k$ is algebraically closed.

\vspace{.5em}

\noindent
{\bf Case 2.\, $k$ arbitrary.}\;  Now, let $k$ be a field of characteristic $0$.

\begin{lemma}\label{invariants and base change}
If $V$ is a rational $G$-representation defined over $k$, then 
\bd
\left(V^G\right)_{\bar{k}} = \left(V_{\bar{k}}\right)^{G_{\bar{k}}}.
\ed
\end{lemma}

By Lemma \ref{invariants and base change}, we have 
$(\pi_*^G\cM)_{\bar{k}} = \pi_*^{G_{\bar{k}}}\cM_{\bar{k}}$, and consequently
\begin{equation}\label{bark of adjoint}
\big(\pi^*\pi_*^G\cM\big)_{\bar{k}}=\pi_{\bar{k}}^*\pi_*^{G_{\bar{k}}}\cM_{\bar{k}}.
\end{equation}
  Since the canonical map \eqref{invariantsgivedescent} for 
$\cM_{\bar{k}}$ is an isomorphism by Lemma \ref{base-changed condition}
and the $\bar{k}$-case of the proof, \eqref{bark of adjoint} implies that
the pullback to $\bar{k}$ of $\pi^*\pi_*^G\cM\rightarrow \cM$ is an
isomorphism.  But now $\spec\bar{k}\rightarrow \spec k$ is faithfully flat,
so $\pi^*\pi_*^G\cM\rightarrow \cM$ is itself an isomorphism.  
By Corollary \ref{canonical adjunction}, this completes the proof.\hfill\qedsymbol

\vspace{.5em}

We note that Condition (2) of Theorem \ref{firstprop} cannot be replaced by 
Condition (3) when $k$ is not algebraically closed:
\begin{example}
Let $X=\spec {\mathbf R}[z,z^{-1}]\cong {\mathbf G}_m$, a variety over
${\mathbf R}$ with $G={\mathbf G}_m$ 
acting in the obvious way.  let $x=(z^2+1)$.  Then $G_x = \{\pm 1\}$.  
Let $\cM = \theo_X$.  The action of $G_x$ on the fiber 
$\cM\otimes\theo_X/\fm_x\cong {\mathbf C}$ is identified with the action
of $\on{Gal}({\mathbf C}/{\mathbf R})$.  In particular, although the fiber
is generated by $G_x$-invariants, it does not consist entirely of 
$G_x$-invariants.
\end{example}

\section{Proof of Theorem \ref{secondtheorem}}\label{proofoftwo}

\subsection{Preliminaries}

\begin{prop}\label{replacementlemma}
Suppose $\bE$
is a $G$-equivariant complex of vector bundles on a noetherian affine scheme $Y$ over an algebraically closed field $k$ of characteristic zero.  
Assume that, for some $n\in{\mathbf Z}$, one is given a $G$-equivariant vector bundle $V$ on $Y$ and a $G$-equivariant homomorphism $V\xrightarrow{f} E_n$ so that
\begin{enumerate}
\item $\im (\phi_n\circ f) = \im (\phi_n)$ and
\item the induced homomorphism 
\bd
\Ker (\phi_n\circ f) \rightarrow H^n(\bE)
\ed
is surjective.
\end{enumerate}
Then there is a $G$-equivariant complex $\bE'$ of vector bundles and a $G$-equivariant quasi-isomorphism $\bE' \xrightarrow{q} \bE$ so that
\begin{enumerate}
\item in degrees $j>n$ one has $E'_j = E_j$, and $q_j: E_j \rightarrow E_j$ is the identity map, and
\item in degree $n$ the quasi-isomorphism $q$ restricts to $V\xrightarrow{f} E_n$.
\end{enumerate}
\end{prop}
\begin{proof}
We produce $\bE'$ recursively, starting from:
\bd
\xymatrix{ & V\ar[d]^{f} \ar[r]^{\phi_n\circ f} & E_{n+1} \ar[d] \ar[r]^{\phi_{n+1}} & E_{n+2} \ar[d] \ar[r] & \dots \\
E_{n-1} \ar[r]^{\phi_{n-1}} & E_n \ar[r]^{\phi_n} & E_{n+1} \ar[r]^{\phi_{n+1}} & E_{n+2} \ar[r] & \dots}
\ed
So, suppose we are given a $G$-equivariant morphism
\bd
\xymatrix{ & E_j'\ar[d]^{q_j} \ar[r]^{\phi_j'} & E_{j+1}' \ar[d] \ar[r]^{\phi_{j+1}'} & E_{j+2} \ar[d] \ar[r] & \dots \\
E_{j-1} \ar[r]^{\phi_{j-1}} & E_j \ar[r]^{\phi_j} & E_{j+1} \ar[r]^{\phi_{j+1}} & E_{j+2} \ar[r] & \dots}
\ed
which induces
\begin{enumerate}
\item[(i)] an isomorphism on $H^i$ for $i\geq j+1$, and
\item[(ii)] a surjection on $H^j$. 
\end{enumerate}
  Let
$\widetilde{E}_{j-1}' = E_{j-1} \underset{E_j}{\times} E_j'$.
  By construction, $\widetilde{E}_{j-1}'$ is equipped with $G$-equivariant
morphisms 
\bd
E_{j-1}\xleftarrow{\widetilde{q}_{j-1}}\widetilde{E}_{j-1}'\xrightarrow{\widetilde{\phi}_{j-1}'} E_j'
\ed
 for which $\phi_{j-1}\circ \widetilde{q}_{j-1} = q_j \circ \widetilde{\phi}_{j-1}'$.  Furthermore, the induced morphism on cohomology in degree $j$ is an isomorphism: by assumption the sheaf $\Ker(\phi_j')$ surjects onto $H^j(\bE)$, and so it is enough to check that the image of $\widetilde{E}_{j-1}'$ in $E_j'$ is the kernel of the map
$\Ker (\phi_j') \rightarrow H^j(\bE)$.
A section $e$ of $\Ker(\phi_j')$ goes to zero in $H^j(\bE)$, however, exactly if there is some $e'$ in $E_{j-1}$ for which  $\phi_{j-1}(e') = q_j(e)$, and then the section $(e',e)$ lies in $\widetilde{E}_{j-1}'$ and maps to $e$.  

In addition, it is easy to see that the kernel of the map
\begin{equation}\label{gettingthekernel}
\widetilde{E}_{j-1}' \rightarrow E_j'
\end{equation}
surjects onto $\Ker(\phi_{j-1})$ (if $e$ is a section of the kernel, then $(e,0)$ lies in the kernel of \eqref{gettingthekernel}) and hence onto $H^{j-1}(\bE)$.  Thus, to complete the proof it will be enough to produce a $G$-equivariant vector bundle $E_{j-1}'$ and a $G$-equivariant surjective morphism of coherent sheaves $E_{j-1}' \rightarrow \widetilde{E}_{j-1}'$.  But the existence of such a vector bundle and morphism is guaranteed by Lemma \ref{canfindsurjection}. 
It follows that our extension of $\bE'$ now satisfies (i) and (ii) with $j$ 
replaced by $j-1$, and continuing our recursion gives the desired complex.
\end{proof}

\begin{lemma}\label{descent degrees}
Suppose that 
\begin{equation}\label{triangle}
\bE'\rightarrow \bE\rightarrow \bE''\xrightarrow{[1]}
\end{equation}
is an exact triangle of $G$-equivariant complexes.  Suppose that the canonical
map \eqref{derived cat canonical} is 
\begin{enumerate}
\item an isomorphism on $H^i$ for $\bE'$ for all $i$, and
\item an isomorphism on $H^i$ for $\bE''$ in degrees $i\geq n$ and a surjection
in degree $i=n-1$.  
\end{enumerate}
Then the canonical map \eqref{derived cat canonical} is
\begin{enumerate}
\item[(a)] an isomorphism on
$H^i$ for $\bE$ in degrees $i\geq n$ and
\item[(b)] a surjection on $H^{n-1}$ for $\bE$.
\end{enumerate}
\end{lemma}
\begin{proof}
This is a restatement of the Five Lemma for 
\eqref{triangle}.
\end{proof}
\begin{lemma}\label{condition}
Given an exact triangle of the form \eqref{triangle}, if two of
$\bE$, $\bE'$, and $\bE''$ satisfy condition (2) of Theorem 
\ref{secondtheorem}, then so does the third.
\end{lemma}
\begin{proof}
This is immediate from the long exact cohomology sequence.
\end{proof}
\begin{lemma}\label{surjective map}
Suppose that $\cM$ is a $G$-equivariant coherent sheaf on $X$, and that 
$\bV$ is a $G$-equivariant vector bundle on $X$ that descends to
$X/\!\!/G$.  Suppose that there is a $G$-equivariant surjective 
map $\bV\rightarrow\cM$.  Then the natural map \eqref{invariantsgivedescent}
for $\cM$ is surjective.
\end{lemma}
\begin{proof}
We have a surjective map
\bd
\pi^*\pi_*^G\bV = \bV \rightarrow \cM
\ed
that factors through $\pi^*\pi_*^G\cM$.
\end{proof}
\begin{cor}\label{last term}
Suppose that $\bM$ is a coherent
 complex concentrated in degrees $(-\infty,m]$ that
satisfies condition (2) of Theorem \ref{secondtheorem}.  Then the
canonical map \eqref{derived cat canonical} is surjective on $H^m$.
\end{cor}
\begin{proof}
By hypothesis and
right-exactness of tensor product, $H^m(\bM)$ satisfies condition 
(2) of Theorem \ref{secondtheorem} as well.  Hence, applying 
Lemma \ref{canfindsurjection} to $H^m(\bM)$, there are a $G$-equivariant
vector bundle $\bV$ that descends to $X/\!\!/G$ and an equivariant
surjection $\bV\rightarrow H^m(\bM)$.  Lemma \ref{surjective map} then
implies that $\pi^*\pi_*^G H^m(\bM)\rightarrow H^m(\bM)$ is surjective.  

We now use the exact triangle
\bd
\tau_{\leq m-1} \bM\rightarrow \bM\rightarrow H^m(\bM)\xrightarrow{[1]}.
\ed
Applying ${\mathbf L}\pi^*\pi_*^G$ to it and taking the long exact cohomology 
sequence, we get an exact sequence
\bd
\dots\rightarrow H^m({\mathbf L}\pi^*\pi_*^G\bM)\rightarrow
\pi^*\pi_*^G H^m(\bM)\rightarrow H^{m-1}({\mathbf L}\pi^*\pi_*^G\tau_{\leq m-1}\bM)\rightarrow\dots.
\ed
But the last group above is zero by degree considerations, so the map
$H^m({\mathbf L}\pi^*\pi_*^G\bM)\rightarrow \pi^*\pi_*^GH^m(\bM)$ is
surjective.  Combining this with the conclusion of the previous paragraph
yields the corollary.
\end{proof}

\subsection{Necessity of the Criterion}

Suppose that $\bE$
is a $G$-equivariant bounded above complex of vector bundles on $X$ that descends to $X/\!\!/G$. 

Let $x\in X$ be a closed point in a closed $G$-orbit.  Then
\bd
\bE\overset{\bL}{\otimes} \theo_X/\fm_x = \big(\bL\pi^*\pi_*^G\bE\big)
\overset{\bL}{\otimes} \theo_X/\fm_x
= \left(\pi_*^G\bE\overset{\bL}{\otimes} \theo_{X/\!\!/G}/\fm_{\pi(x)}\right)
\overset{\bL}{\otimes}\theo_X/\fm_x.
\ed
Now $\theo_{X/\!\!/G}/\fm_{\pi(x)}\rightarrow \theo_X/\fm_x$ is flat, so
\bd
H^i(\bE\overset{\bL}{\otimes}\theo_X/\fm_x) = H^i(\pi_*^G\bE\overset{\bL}{\otimes} \theo_{X/\!\!/G}/\fm_{\pi(x)})\otimes\theo_X/\fm_x.
\ed
The right-hand side is generated over $\theo_X/\fm_x$ by 
 the $G_x$-invariants $H^i(\pi_*^G\bE\overset{\bL}{\otimes}\theo_{X/\!\!/G}/\fm_{\pi(x)})$.

\subsection{Sufficiency of the Criterion}
\mbox{}

\noindent
{\bf Case 1.\,  $k$ algebraically closed.}\;
We prove the sufficiency of Condition (2) by an induction.  Consider the following statement:

\vspace{.7em}

\hangup\noindent
$P(k)$: Suppose that $\bE$ is a $G$-equivariant complex that satisfies 
Condition (2) of Theorem \ref{secondtheorem} and that $\bE$ is concentrated in
degrees $(-\infty, m]$.  Then \eqref{derived cat canonical} induces an
isomorphism on $H^i$ for $i\geq m-k$ and a surjection on $H^{m-k-1}$.

\vspace{.7em}

To complete the proof of the theorem, it suffices to prove $P(k)$ for all $k$.
Moreover, $P(-2)$ holds trivially.

 Suppose, by way of inductive hypothesis,
 that $P(k)$ holds (that is, for all $\bE$).  Let
\bd
{\mathbf E}: \hspace{1.5em} \dots\rightarrow E_n \xrightarrow{\phi_n} E_{n+1} \rightarrow \dots \rightarrow E_{-2} \xrightarrow{\phi_{-2}} E_{-1} 
\xrightarrow{\phi_{-1}} E_0 \rightarrow 0
\ed
be a $G$-equivariant complex of vector bundles on $X$.
By Corollary \ref{last term} and
Lemma \ref{canfindsurjection}, there is a $G$-equivariant vector bundle
$\bE_0'$ that descends to $X/\!\!/G$  and a surjective homomorphism
$\bE_0'\rightarrow H^0(\bE)$; in light of Remark \ref{reduce to affine},
we may assume that this homomorphism lifts to an equivariant homomorphism
$\bE_0'\rightarrow \bE$.

  Now, applying Proposition \ref{replacementlemma},
we find that there is a complex $\bE'$ that is $G$-equivariantly 
quasi-isomorphic to $\bE$, such that $\bE'$ is concentrated in degrees
$(-\infty, 0]$ and $\bE'_0$ descends.  We get an exact triangle
\begin{equation}\label{E' triangle}
\bE'_0\rightarrow \bE'\rightarrow \sigma_{\leq -1}\bE'\xrightarrow{[1]}
\end{equation}
for which $\bE'_0$ descends to $X/\!\!/G$ and $\sigma_{\leq -1}\bE'$ is
concentrated in degrees $(-\infty, -1]$.  By Lemma \ref{condition}, all three
terms in \eqref{E' triangle} satisfy condition (2) of the theorem.  

We wish to apply Lemma \ref{descent degrees} to 
\eqref{E' triangle}.
  Since $\bE'_0$ descends, condition (1) of Lemma \ref{descent 
degrees} is satisfied.  The inductive hypothesis tells us that
 condition (2) of Lemma \ref{descent degrees}, applied to $\sigma_{\leq -1}\bE'$, holds for $n=-k-1$.   Thus, Lemma \ref{descent degrees} implies that
$P(k+1)$ holds 
for $\bE'$.  Consequently, $P(k+1)$ holds for $\bE$.  Since this conclusion
holds for arbitrary $\bE$,  this completes the 
inductive step, and thus the proof of the theorem when $k$ is algebraically
closed.

\vspace{.5em}

\noindent
{\bf Case 2.\, $k$ arbitrary.}\;
Consider the morphism \eqref{derived cat canonical}.  Tensoring with 
$\bar{k}$ and using \eqref{bark of adjoint}, we obtain the morphism
\bd
(\bL\pi^*\pi_*^G\bE)_{\bar{k}} = \bL\pi^*\pi_*^{G_{\bar{k}}}\bE_{\bar{k}}
\rightarrow \bE_{\bar{k}}.
\ed
By Lemma \ref{base-changed condition} and the algebraically closed case, 
this is a quasi-isomorphism.  But $\spec\bar{k}\rightarrow\spec k$ is 
faithfully flat, so \eqref{derived cat canonical} is also a quasi-isomorphism
for $\bE$.  By Corollary \ref{canonical adjunction}, this completes the
proof.\hfill\qedsymbol

\section{Some Examples}\label{sympl section}
In this section, we discuss two examples that have motivated the author's interest in the descent questions considered here.
\subsection{Cotangent Complex in the Smooth Setting}
Let $X$ be a smooth complex variety with an action of a reductive group $G$ and good quotient $\pi: X\rightarrow X/\!\!/G$.  The equivariant cotangent 
complex of $X$ is the complex 
\bd
{\mathbb L}: \hspace{1em}\Omega^1_X \rightarrow \theo_X\otimes {\mathfrak g}^*
\ed
 in which, on fibers, the map 
is given by the dual of the infinitesimal action
of ${\mathfrak g} = \on{Lie}(G)$ on $X$.  This complex is the pullback of the cotangent complex of the stack quotient $[X/G]$ along
the projection $X\rightarrow [X/G]$, hence it measures ``singularities in the equivariant geometry of $X$.''

 Suppose that ${\mathbb L}$ descends to $X/\!\!/G$ on a neighborhood of the image 
$\pi(x)$ of a closed point $x\in X$ lying in a closed $G$-orbit.  Then the stabilizer $G_x$ must act trivially on the cohomologies
$H^i({\mathbb L}_x)$ of the complex $\Omega^1_{X,x}\rightarrow {\mathfrak g}^*$, where $\Omega^1_{X,x}$ denotes the cotangent space to $X$ at
$x$.  Now $H^1({\mathbb L}_x) = {\mathfrak g}_x^*$ is the dual of the Lie algebra of $G_x$ (with the coadjoint action of $G_x$) and
$H^0({\mathbb L}_x) = (T_{X,x}/{\mathfrak g})^\vee = N^\vee_{G\cdot x/X}(x)$ is the conormal space to the orbit $G\cdot x$ at $x$.  If $G_x$ 
acts trivially on $N^\vee_{G\cdot x/X}(x)$, it follows from
Luna's \'Etale Slice Theorem that there then exists a smooth slice $S$ through $x$ for the $G$-action on $X$ on which $G_x$
acts trivially.  Moreover, since $X$ is \'etale-locally isomorphic to $G\times_{G_x} S$ in a neighborhood of $x$, we may conclude
that $X\cong (G/G_x)\times S$ \'etale-locally
and $G$-equivariantly near $x$.  In particular, $x$ lies in the open stratum $U$ of $X$ on which 
the orbit type is constant, i.e. every $p\in U$ has stabilizer conjugate to $G_x$.  Summarizing:
\begin{prop}
Suppose $X$ is a smooth complex variety with an action of a reductive group $G$ with good quotient
$\pi: X\rightarrow X/\!\!/G$.  Let $U\subset X$ denote the open stratum described above.  Then $\pi(U)$ is the largest open subset of the
good quotient $X/\!\!/G$ 
to which the equivariant cotangent complex ${\mathbb L}$ descends.
\end{prop}

\subsection{Cotangent Complex in the Symplectic Setting}
As we mentioned in the introduction, Theorem \ref{secondtheorem} was originally motivated by the following picture coming from
(algebraic) symplectic geometry.  Let $M$ be a smooth, affine complex variety with an action of a reductive group $G$.  Let 
$\mu: T^*M \rightarrow {\mathfrak g}^*$ be a moment map for the induced action on the cotangent bundle.  Suppose, for simplicity, 
that $\mu$ makes $N = \mu^{-1}(0)$ a complete intersection, so that the quotient stack ${\mathcal X} = [N/G]$ is a complete intersection and its
cotangent complex ${\mathbb L}_{\mathcal X}$ is concentrated in $[-1,1]$.
Moreover, ${\mathbb L}_{\mathcal X}$
comes equipped with a map $\omega: {\mathbb L}_{\mathcal X} \rightarrow {\mathbb L}_{\mathcal X}^\vee$ that is antisymmetric and nondegenerate (that is, it pairs $H^0$ 
nondegenerately with itself and $H^1$ nondegenerately with $H^{-1}$---here antisymmetry is taken in the graded sense and the 
nondegeneracy is a statement about the pairing on cohomology of fibers).  
The quotient stack $[N/G]$ equipped with this structure---an antisymmetric, closed, and nondegenerate pairing on its cotangent complex---is what should be properly considered a (local complete intersection) symplectic stack.\footnote{If a stack $M$ is not lci, 
then the cotangent complex will lie in degrees $(-\infty, 1]$ but {\em not} in $[-1,1]$.  In that case, it is reasonable to hope that 
${\mathbb L}_{\mathcal X}$ will admit a nondegenerate pairing for a 
suitable derived enhancement $\widetilde{\mathcal X}$ of ${\mathcal X}$.}

It is natural to wonder whether one can extract information about the singularities of the quotient {\em space} $X = N/\!\!/G$ 
from ${\mathbb L}_{\mathcal X}$: for example, the dimensions of cohomologies of ${\mathbb L}_{\mathcal X}$ stratify ${\mathcal X}$
(and hence $N$) by dimensions of stabilizers, which induces a stratification of $X$ as well.    In this light, it would be optimal if
the pullback of ${\mathbb L}_{\mathcal X}$ to $N$, i.e. the $G$-equivariant cotangent complex of $N$, descended to $X$.  To check 
this, one can use the following description of the cotangent complex.
First, the cotangent complex
of $N$ is given by the pullback map on $1$-forms,
\bd
\theo_N\otimes {\mathfrak g}\xrightarrow{\mu^*} (\Omega^1_{T^*M})|_N.
\ed
The cotangent complex of the stack $[N/G]$ is then obtained by descending the $G$-equivariant complex on $N$,
\bd
\widetilde{T}^*: \hspace{2em}\theo_N\otimes {\mathfrak g} \xrightarrow{\mu^*} (\Omega^1_{T^*M})|_N \xrightarrow{da^*} \theo_N\otimes {\mathfrak g}^*,
\ed
where $da^*$ is the dual of the infinitesimal action map $da: \theo_N\otimes {\mathfrak g}\rightarrow T_{T^*M}$.

Unfortunately, examples indicate that this $G$-equivariant perfect complex on $N$ probably rarely (or maybe never) descends to
$N/\!\!/G$; but this failure of descent seems to be an interesting phenomenon worthy of further consideration.  The closures of
minimal nilpotent orbits in semisimple groups provide an interesting class of examples.  In type $A$, such orbit closures are symplectic
reductions of $T^*{\mathbb C}^{n+1}$ by the ${\mathbb C}^*$-action induced from scaling on ${\mathbb C}^{n+1}$.  The quotient
is a subvariety of ${\mathfrak sl}_{n+1}$, namely the set of nilpotent matrices of rank less than or equal to $1$.  Indeed, 
the moment map 
\bd
\mu: T^*{\mathbb C}^{n+1} = {\mathbb C}^{n+1}\times ({\mathbb C}^{n+1})^*\rightarrow {\mathbb C}
\ed
is $\mu(i,j) = j(i)$, which is equivalent to $\on{trace}(ij) = 0$; the map $\mu^{-1}(0)\rightarrow {\mathfrak sl}_{n+1}$ takes
$(i,j)\mapsto ij$.  The orbit closure $\overline{\mathcal O} = \mu^{-1}(0)/\!\!/{\mathbb C}^*$ has an isolated singularity at the origin in ${\mathfrak sl}_{n+1}$, 
and it is easy to see using Theorem \ref{secondtheorem} that the equivariant cotangent complex of $\mu^{-1}(0)$ does not descend
to $\overline{\mathcal O}$
at the singularity.  However, $\overline{\mathcal O}$ has a symplectic resolution of singularities by $T^*{\mathbb P}^n$ 
\cite{MR2003k:14061}, realized as the quotient of the subset
\bd
\mu^{-1}(0)^s = \mu^{-1}(0)\cap \left(({\mathbb C}^{n+1}\smallsetminus \{0\})\times ({\mathbb C}^{n+1})^*\right)
\ed
by the action of ${\mathbb C}^*$.  Every point of $\mu^{-1}(0)^s$ has reductive stabilizer---in fact, the action is free---and so
the restriction of the equivariant cotangent complex naturally descends to $T^*{\mathbb P}^n = \mu^{-1}(0)^s/{\mathbb C}^*$.  

This is, of course, a rather trivial example, but it is intriguing to compare it to cases in which a symplectic resolution is known
{\em not} to exist: for example, minimal nilpotent orbit closures in type $C$ \cite{MR2003k:14061}.  In this case, the orbit
closure is a quotient of ${\mathbb C}^{2n}$ by the action of ${\mathbb Z}/2{\mathbb Z}$, and again 
the cotangent complex---in this case, just the cotangent bundle---does not descend.  However, in this example there is
no clear method to obtain a quotient to which the cotangent complex does descend: the largest open set on which this does happen,
which again consists of free orbits, does not map surjectively to 
${\mathbb C}^{2n}/\!\!/({\mathbb Z}/2{\mathbb Z})$, and so does not provide a symplectic resolution of the minimal orbit closure.

It would be interesting to know whether one can explain the existence of the symplectic resolution in terms of the behavior of the 
cotangent complex.  It is conjectured (cf. \cite{MR2003k:14061}) that every birational contraction of a smooth symplectic variety is locally modelled on a symplectic resolution
of some nilpotent orbit closure.  Thus, one may hope that understanding the behavior of the cotangent complex for nilpotent orbit
closures may shed light on existence of symplectic resolutions.

\bibliographystyle{alpha}
\def\cprime{$'$} \def\cprime{$'$} \def\cprime{$'$} \def\cprime{$'$}
  \def\cprime{$'$} \def\cprime{$'$} \def\cprime{$'$} \def\cprime{$'$}
  \def\cprime{$'$} \def\cprime{$'$} \def\cprime{$'$} \def\cprime{$'$}
  \def\cprime{$'$} \def\cprime{$'$} \def\cprime{$'$} \def\cprime{$'$}
  \def\cprime{$'$} \def\cprime{$'$} \def\cprime{$'$} \def\cprime{$'$}
  \def\cprime{$'$} \def\cprime{$'$} \def\cprime{$'$} \def\cprime{$'$}
  \def\cprime{$'$} \def\cprime{$'$} \def\cprime{$'$} \def\cprime{$'$}
  \def\cprime{$'$} \def\cprime{$'$} \def\cprime{$'$} \def\cprime{$'$}
  \def\cprime{$'$} \def\cprime{$'$} \def\cprime{$'$}

\end{document}